\documentclass[11pt,a4paper,leqno]{amsart}
\usepackage[latin1]{inputenc}
\usepackage[T1]{fontenc}
\usepackage{amsfonts}
\usepackage{amsmath}
\usepackage{amssymb}
\usepackage{eurosym}
\usepackage{mathrsfs}
\usepackage{palatino}
\usepackage{color}
\usepackage{esint}
\usepackage{url}
\usepackage{verbatim}
\usepackage{graphicx}
\usepackage{enumerate}
\usepackage[pagebackref,hypertexnames=false, colorlinks, citecolor=blue, linkcolor=blue, urlcolor=red]{hyperref}

\allowdisplaybreaks[4]

\numberwithin{equation}{section}

\swapnumbers
\theoremstyle{plain}
\newtheorem{thm}[equation]{Theorem}
\newtheorem{lem}[equation]{Lemma}
\newtheorem{prop}[equation]{Proposition}

\theoremstyle{definition}

\theoremstyle{remark}
\newtheorem{rem}[equation]{Remark}

 \email{zynn@tju.edu.cn}
 \email{fandashan2@zjnu.edu.cn}
\makeatletter
\@namedef{subjclassname@2020}{  \textup{2020} Mathematics Subject Classification}
\makeatother
\newcommand{\func}[1]{\operatorname{#1}}

\begin{document}
\title{Variational Estimates for the Bilinear Hilbert Transform}
\author{Yunan Zeng}
\author{Dashan Fan}
\address[Y.Z.]{Center for Applied Mathematics, Tianjin University, Weijin
Road 92, 300072 Tianjin, China}
\address[D.F.]{Department of Mathematical Science, University of
Wisconsin-Milwaukee, Milwaukee, WI 53201, USA}

\begin{abstract}
    We prove that for every $q>2$, the $q$-variation operator associated with the bilinear Hilbert transform
is bounded from $L^r(\mathbb R)\times L^s(\mathbb R)$ to
$L^p(\mathbb R)$ whenever $1<r,s<\infty$, $p>2/3$, and
$1/p=1/r+1/s$. 
\end{abstract}

\maketitle

\allowdisplaybreaks


\section{Introduction}
L\'epingle \cite{Lepingle1976} first established variational inequalities for martingales. Bourgain later used variational methods
to study the pointwise convergence of ergodic averages along arithmetic sets \cite{Bourgain1989}. In recent years, variation operators have received considerable
attention in harmonic analysis and ergodic theory. Subsequent work established variational estimates for linear averaging
operators \cite{JSW2008,JRW2003,MST2017,BMSW2018}, truncated singular integrals \cite{CJRW2000,CJRW2003,DHL2017,CH2023}, bilinear paraproducts
\cite{DMT}, and bilinear averaging operators \cite{Liu2020,DHW2019}.

In this paper, we study variational estimates for the bilinear Hilbert
transform, which is defined by 
\begin{equation*}
H(f,g)(x)=p.v.\int_{\mathbb{R}}f(x-y)g(x+y)\frac{dy}{y}.
\end{equation*}%
Lacey and Thiele \cite{LT1997,LT1999} established the classical
strong-type bound
\begin{equation*}
    \|H(f,g)\|_{L^p(\mathbb R)}\lesssim_{r,s,p}\|f\|_{L^r(\mathbb R)}\|g\|_{L^s(\mathbb R)}
\end{equation*}
for $1<r,s\leq\infty$ and $2/3<p<\infty$ with $1/p=1/r+1/s$.
More recently, Di Plinio \cite{DiPlinio2026} established the
weak-type endpoint estimate
\begin{equation*}
\|H(f,g)\|_{L^{2/3,\infty}(\mathbb R)}\lesssim_{r,s}\|f\|_{L^r(\mathbb R)}\|g\|_{L^s(\mathbb R)}
\end{equation*}
for $1<r,s<2$ with $1/r+1/s=3/2$. Hu \cite{Hu2026} also gave an alternative proof.

For $\varepsilon >0$, we define the truncation of $H(f,g)$ as 
\begin{equation*}
H_{\varepsilon }(f,g)(x)=\int_{|y|>\varepsilon }f(x-y)g(x+y)\,\frac{dy}{y}.
\end{equation*}%
The associated maximal truncation is
\begin{equation*}
   H^*(f,g)(x)=\sup_{\varepsilon>0}|H_\varepsilon(f,g)(x)|. 
\end{equation*}
Lacey's maximal truncation theorem \cite{Lacey} implies that
\begin{equation*}
    \|H^*(f,g)\|_{L^p(\mathbb R)}\lesssim_{r,s,p}\|f\|_{L^r(\mathbb R)}\|g\|_{L^s(\mathbb R)}
\end{equation*}
whenever $1<r,s\leq\infty$, $2/3<p<\infty$, and $1/p=1/r+1/s$. 

Demeter \cite{Demeter2007} used an oscillation estimate to prove
almost-everywhere convergence for the ergodic bilinear Hilbert
transform. Berkson and Demeter \cite{BerksonDemeter2010} studied
the pointwise convergence of its spatial truncations.
Do, Muscalu, and Thiele \cite{DMT2017} established variational
estimates for the bilinear iterated Fourier integral with
frequency truncations.

For $1\le q<\infty$, we define the $q$-variation operator
associated with the spatial truncations $\{H_\varepsilon\}_{\varepsilon>0}$ by
\begin{equation*}
V_{q}\big(H_{\varepsilon }(f,g):\varepsilon >0\big)(x):=\sup_{0<\varepsilon
_{0}<\ldots <\varepsilon _{N}}\Big(\sum_{k=0}^{N-1}\big|H_{\varepsilon
_{k+1}}(f,g)(x)-H_{\varepsilon _{k}}(f,g)(x)\big|^{q}\Big)^{\frac{1}{q}}.
\end{equation*}%
Our main result is the following variational estimate, which also
yields $L^r\times L^s\to L^p$ bounds for both $H$ and its maximal
truncation $H^*$.
\begin{thm}
Let $1<r,s<\infty,p>{2}/{3}$ and ${1}/{p}=1/s+1/r$. Then, for every $q>2$ we
have 
\begin{equation*}
\Big\|V_q\big(H_{\varepsilon}(f,g):\varepsilon>0\big)\Big\|_{L^p(\mathbb{R}%
)}\leq C_{p,q,r,s}\|f\|_{L^r(\mathbb{R})}\|g\|_{L^s(\mathbb{R})}.
\end{equation*}
\end{thm}

By density, it suffices to prove the estimate for $f,g\in \mathcal{S}(%
\mathbb{R})$. We begin the proof by decomposing the operator \ $%
H_{\varepsilon }$. Choose an even function $\chi \in C^{\infty }(\mathbb{R})$
such that $0\leq \chi \leq 1$, $\chi $ is nondecreasing on $[0,\infty )$,
and 
\begin{equation*}
\chi (u)=0\quad \text{for }|u|\leq 1,\qquad \chi (u)=1\quad \text{for }%
|u|\geq 2.
\end{equation*}%
Define the smooth truncation 
\begin{equation*}
H_{\varepsilon }^{\chi }(f,g)(x):=\int_{\mathbb{R}}f(x-y)g(x+y)\frac{\chi
(y/\varepsilon )}{y}dy
\end{equation*}%
and the error term 
\begin{equation*}
E_{\varepsilon }(f,g)(x):=\int_{\mathbb{R}}f(x-y)g(x+y)\frac{%
1_{\{|y|>\varepsilon \}}-\chi (y/\varepsilon )}{y}dy.
\end{equation*}%
Then 
\begin{equation*}
H_{\varepsilon }(f,g)(x)=H_{\varepsilon }^{\chi }(f,g)(x)+E_{\varepsilon
}(f,g)(x).
\end{equation*}%
We will establish variational estimates for these two terms. The theorem
then follows from the pointwise triangle inequality for $V_{q}$ and, when $%
p<1$, the $p$-subadditivity of the $L^{p}$ quasi-norm. 

For the smooth truncations, we verify the multiplier conditions in
\cite[Theorem~1.3]{DOP} and apply the resulting variational estimate.
The main difficulty is to control the error term, which reduces to
variational estimates for bilinear averages. Do, Oberlin, and
Palsson \cite[Theorem~2.1]{DOP} obtained such estimates for
sufficiently large $q$. Becker and Durcik \cite[Theorems~1.7 and~1.8]{BD}
proved the long-variation estimate for every $q>2$ and a short-variation
estimate with an outer $\ell^2$ norm. When $2/3<p<1$, their short-variation
estimate requires $1/q<3/2-1/p$. We use a probabilistic argument to
establish the short-variation estimate with an outer $\ell^q$ norm
for every $q>2$ in the range $p>2/3$. The long and short variations and the bilinear averaging operators
will be defined in Section~3.

This paper is organized as follows: Section 2 establishes variational estimates for the smooth truncations. 
Section 3 proves variational estimates for bilinear averages. 
In Section 4, we represent the error kernel as an integral against a finite signed measure and apply the estimates from Section 3.

\section{Variational estimates for smooth truncation}

\label{set 2} \ In this section we prove variational estimates for the
smooth-truncation operator $H_{\varepsilon }^{\chi }$. This operator belongs
to a family of bilinear operators whose variational theory was developed by
Do, Oberlin and Palsson. We verify the required Fourier decay conditions for
the kernel $K(u)=\chi (u)/u$, which yields the $L^{p}$ bound for $V_{q}\big(%
H_{\varepsilon }^{\chi }(f,g):\varepsilon >0\big)$ by applying their
theorem. 

\begin{thm}
\label{thm smooth} Let $1<r,s<\infty $, $2/3<p<\infty $ and $1/p=1/r+1/s$.
Then for every $q>2$, we have 
\begin{equation*}
\Big\|V_{q}\big(H_{\varepsilon }^{\chi }(f,g):\varepsilon >0\big)\Big\|%
_{L^{p}(\mathbb{R})}\leq C_{p,q,r,s}\Vert f\Vert _{L^{r}(\mathbb{R})}\Vert
g\Vert _{L^{s}(\mathbb{R})}.
\end{equation*}
\end{thm}

To prove the above theorem, we recall the following $1$-dimension version of
the variational inequality by Do, Oberlin, and Palsson \cite[Theorem~1.3]%
{DOP}.

\begin{thm}[{Do--Oberlin--Palsson {\protect\cite[Theorem~1.3]{DOP}}}]
\label{thm Do} Let $1<r,s<\infty$, $2/3<p<\infty$, and $2<q<\infty$ satisfy $%
1/p=1/r+1/s.$ There exists an integer $N_0=N_0(p,q,r,s)\ge1$ with the
following property.

Suppose that $K:\mathbb{R}\rightarrow \mathbb{C}$ satisfies the following
estimates uniformly for $\xi \neq 0$: 
\begin{equation*}
|\widehat{K}(\xi )|\leq A_{0}\min \{1,|\xi |^{-1}\},
\end{equation*}%
and 
\begin{equation*}
|\frac{d^{m}}{d\xi ^{m}}\widehat{K}(\xi )|\leq A_{m}\min \{|\xi |^{1-m},|\xi
|^{-m-1}\},\qquad 1\leq m\leq N_{0}.
\end{equation*}%
Set 
\begin{equation*}
B_{t}(f,g)(x):=\frac{1}{t}\int_{\mathbb{R}}f(x-y)g(x+y)K(y/t)\,dy.
\end{equation*}%
Then 
\begin{equation*}
\left\Vert V_{q}\bigl(B_{t}(f,g):t>0\bigr)\right\Vert _{L^{p}(\mathbb{R}%
)}\leq C\Vert f\Vert _{L^{r}(\mathbb{R})}\Vert g\Vert _{L^{s}(\mathbb{R})},
\end{equation*}%
where $C$ depends only on $p,q,r,s,N_{0}$ and $A_{0},\ldots ,A_{N_{0}}$.
\end{thm}

In view of this theorem, it suffices to verify the above conditions for the
kernel $K(u)=\chi(u)/u$ associated with the smooth truncations.

\begin{proof}[Proof of Theorem~\protect\ref{thm smooth}]
First, it is trivial to see that $|\widehat{K}(\xi )|\lesssim 1$.

Next, we prove $|\widehat{K}(\xi )|\lesssim |\xi |^{-1}$. For $|u|>2$, we
have 
\begin{equation*}
K^{(N)}(u)=\frac{(-1)^{N}N!}{u^{N+1}}.
\end{equation*}%
and $K^{(N)}\in L^{1}(\mathbb{R})$ for every integer $N\geq 1$. Repeating
integration by parts yields 
\begin{equation*}
\int_{-\infty }^{\infty }K(u)e^{-iu\xi }du=\frac{1}{(i\xi )^{N}}%
\int_{-\infty }^{\infty }K^{(N)}(u)e^{-iu\xi }du,
\end{equation*}%
for any $N\geq 1.$

It remains to show that 
\begin{equation*}
\left\vert \frac{d^{m}}{d\xi ^{m}}\widehat{K}(\xi )\right\vert \lesssim
_{m}\min \left\{ |\xi |^{1-m},|\xi |^{-m-1}\right\} ,\qquad \xi \neq 0,\quad
1\leq m\leq N_{0}.
\end{equation*}%
\ Noting 
\begin{equation*}
K\left( u\right) =p.v.\left( \frac{1}{u}\right) -\frac{1-\chi (u)}{u}
\end{equation*}%
we have that 
\begin{equation*}
\widehat{K}(\xi )=-i\pi \text{ }sgn(\xi )-\widehat{R}\left( \xi \right) 
\end{equation*}%
where 
\begin{equation*}
\widehat{R}\left( \xi \right) =2i\int_{0}^{\infty }\left( 1-\chi (u)\right) 
\frac{\sin u\xi }{u}du.
\end{equation*}%
Since 
\begin{equation*}
\frac{d^{m}}{d\xi ^{m}}sgn(\xi )=0
\end{equation*}%
for any $\xi \neq 0,$ it suffices to check $\frac{d^{m}}{d\xi ^{m}}\widehat{R%
}\left( \xi \right) .$

If $0<\left\vert \xi \right\vert \leq 1$, we have that 
\begin{equation*}
\frac{d^{m}}{d\xi ^{m}}\widehat{R}\left( \xi \right) =2i\int_{0}^{\infty
}\left( 1-\chi (u)\right) u^{m-1}\sin ^{\left( m\right) }u\xi \text{ }du,
\end{equation*}%
which yields that when $0<\left\vert \xi \right\vert \leq 1,$ 
\begin{equation*}
\left\vert \frac{d^{m}}{d\xi ^{m}}\widehat{R}\left( \xi \right) \right\vert
\preceq 2\int_{0}^{2}u^{m-1}du\preceq 1.
\end{equation*}%
Next, we estimate the case $\left\vert \xi \right\vert >1.$ Without loss of
generality, we assume $\xi >1,$ since the case of $\xi <-1\ $can be treated
similarly. By changing variables, we write 
\begin{eqnarray*}
\widehat{R}\left( \xi \right)  &=&2i\int_{0}^{\infty }\left( 1-\chi (\frac{u%
}{\xi })\right) \frac{\sin u}{u}du \\
&=&2i\int_{0}^{1}\frac{\sin u}{u}du+2i\text{ }J\left( \xi \right) ,
\end{eqnarray*}%
where 
\begin{equation*}
J\left( \xi \right) =\int_{1}^{\infty }\left( 1-\chi (\frac{u}{\xi })\right) 
\frac{\sin u}{u}du.
\end{equation*}%
It remains to estimate $\frac{d^{m}}{d\xi ^{m}}J\left( \xi \right) .$\ An
integration by parts yields 
\begin{equation*}
J\left( \xi \right) =\cos 1-J_{1}\left( \xi \right) -J_{2}\left( \xi \right)
,
\end{equation*}%
where 
\begin{equation*}
J_{1}\left( \xi \right) =\frac{1}{\xi }\int_{1}^{\infty }\chi ^{\prime }(%
\frac{u}{\xi })\frac{\cos u}{u}du
\end{equation*}%
and 
\begin{equation*}
J_{2}\left( \xi \right) =\int_{1}^{\infty }\left( 1-\chi (\frac{u}{\xi }%
)\right) \frac{\cos u}{u^{2}}du.
\end{equation*}%
Now 
\begin{eqnarray*}
\left\vert \frac{d^{m}}{d\xi ^{m}}J_{2}\left( \xi \right) \right\vert 
&=&\left\vert \int_{1}^{\infty }\partial _{\xi }^{m}\left( 1-\chi (\frac{u}{%
\xi })\right) \frac{\cos u}{u^{2}}du\right\vert  \\
&\lesssim &\frac{1}{\xi ^{m}}\int_{\xi }^{2\xi }\frac{du}{u^{2}}\lesssim 
\frac{1}{\xi ^{m+1}}. \\
&&
\end{eqnarray*}%
Changing variable and observing the support condition of $\chi ^{\prime }$
we now have 
\begin{equation*}
J_{1}\left( \xi \right) =\frac{1}{\xi }\int_{1}^{2}\chi ^{\prime }(u)\frac{%
\cos \xi u}{u}du.
\end{equation*}%
Since $\int_{1}^{2}\chi ^{\prime }(u)\frac{\cos \xi u}{u}du$ is a Schwartz
function, clearly 
\begin{equation*}
\left\vert \frac{d^{m}}{d\xi ^{m}}J_{1}\left( \xi \right) \right\vert
\lesssim \frac{1}{\xi ^{m+1}}.
\end{equation*}

Thus, all the required kernel conditions are satisfied, and the desired
estimate follows from Theorem \ref{thm Do}.
\end{proof}

\section{Variational Estimates for the bilinear averages}

\label{set 4} This section is devoted to variational estimates for the
bilinear averaging operators $A_{t}^{+}$. These averages appear naturally in
the measure-theoretic representation of the error term $E_{\varepsilon }$
from Section 4. Following a standard long-short variation decomposition, the
full q-variation is split into long-range and short-range components. The
long-variation bound is quoted from existing work in the literature, while
the short-variation estimate constitutes the main technical contribution of
this section. To handle the short-variation part, we combine maximal and
square-function estimates with a probabilistic counting stopping-time lemma
for arrays of bilinear averages. Define 
\begin{equation*}
A_{t}^{+}(f,g)(x):=\frac{1}{t}\int_{0}^{t}f(x-y)g(x+y)dy,\quad
A_{t}^{-}(f,g)(x):=\frac{1}{t}\int_{0}^{t}f(x+y)g(x-y)dy.
\end{equation*}%
Recall the long and short $q$-variations: 
\begin{equation*}
V^{L,q}\bigl(A_{t}^{\pm }(f,g)\bigr)(x):=\sup_{j_{0}<\cdots <j_{N},j_{k}\in 
\mathbb{Z}}\Bigl(\sum_{k=0}^{N-1}\big|A_{2^{j_{k+1}}}^{\pm
}(f,g)(x)-A_{2^{j_{k}}}^{\pm }(f,g)(x)\big|^{q}\Bigr)^{\frac{1}{q}},
\end{equation*}%
\begin{equation*}
V^{S,q}\bigl(A_{t}^{\pm }(f,g)\bigr)(x):=\Bigl(\sum_{j\in \mathbb{Z}%
}\sup_{2^{j}\leq t_{0}<\cdots <t_{N}\leq 2^{j+1}}\sum_{k=0}^{N-1}\Big|%
A_{t_{k+1}}^{\pm }(f,g)(x)-A_{t_{k}}^{\pm }(f,g)(x)\Big|^{q}\Bigr)^{\frac{1}{%
q}}.
\end{equation*}%
By the standard long-short variation decomposition, we have 
\begin{equation*}
V_{q}\big(A_{t}^{\pm }(f,g)\big)(x)\lesssim _{q}V^{L,q}\bigl(A_{t}^{\pm
}(f,g)\bigr)(x)+V^{S,q}\bigl(A_{t}^{\pm }(f,g)\bigr)(x).
\end{equation*}

\begin{thm}
\label{Becher} Let $1<r,s<\infty$, $p>2/3$, and $1/p=1/r+1/s$. Then, for $%
q>2 $ we have 
\begin{equation*}
\left\|V_q\big(A_t^{\pm}(f,g):t>0\big)\right\|_{L^p(\mathbb{R})} \le
C_{p,q,r,s}\|f\|_{L^r(\mathbb{R})}\|g\|_{L^s(\mathbb{R})}.
\end{equation*}
\end{thm}

\begin{proof}
Since $A_{t}^{-}(f,g)=A_{t}^{+}(g,f)$, it suffices to consider $%
A_{t}^{+}(f,g)$, and establish the following two propositions.
\end{proof}

\begin{prop}[{\protect\cite[Theorem~1.7]{BD}}]
\label{Becher-long} Let $1<r,s<\infty$, $p>2/3$, and $1/p=1/r+1/s$. Then,
for every $q>2$, we have 
\begin{equation*}
\left\|V^{L,q}\bigl(A_t^\pm(f,g)\bigr)\right\|_{L^p(\mathbb{R})} \le
C_{p,q,r,s}\|f\|_{L^r(\mathbb{R})}\|g\|_{L^s(\mathbb{R})}.
\end{equation*}
\end{prop}

\begin{prop}
\label{Becher-short} Let $1<r,s<\infty$, $p>2/3$, and $1/p=1/r+1/s$. Then,
for every $q>2$ we have 
\begin{equation*}
\left\|V^{S,q}\bigl(A_t^\pm(f,g)\bigr)\right\|_{L^p(\mathbb{R})}\le
C_{p,q,r,s} \|f\|_{L^r(\mathbb{R})}\|g\|_{L^s(\mathbb{R})}.
\end{equation*}
\end{prop}

\begin{rem}
Becker and Durcik \cite[Theorem~1.8]{BD} state the estimate 
\begin{equation*}
\bigg\|\Bigl(\sum_{j\in \mathbb{Z}}V_{q}\bigl(A_{t}^{\pm }(f,g):2^{j}\leq
t\leq 2^{j+1}\bigr)^{2}\Bigr)^{1/2}\bigg\|_{L^{p}(\mathbb{R})}\leq
C_{p,q,r,s}\Vert f\Vert _{L^{r}(\mathbb{R})}\Vert g\Vert _{L^{s}(\mathbb{R})}
\end{equation*}%
under the additional condition $q>\max \left\{ 1,\frac{2p}{3p-2}\right\} $.
For $q>2$, their outer $\ell ^{2}$ norm dominates the outer $\ell ^{q}$ norm
in $V^{S,q}$. Thus, their estimate gives stronger control on the common
range of exponents. When $2/3<p<1$, their condition requires $q>2p/(3p-2)$,
while Proposition~\ref{Becher-short} gives the short-variation estimate with
outer $\ell ^{q}$ for every $q>2$.
\end{rem}

To prove Proposition~\ref{Becher-short}, we decompose the averaging kernel $%
\mathbf{1}_{[0,1]}$ into a smooth term and two families of dyadic pieces
supported near $u=0$ and $u=1$, respectively. The main step is to bound the
short variation of the operators associated with the pieces near $u=1$. The
kernel decomposition is adapted from Becker and Durcik \cite[Lemma~4.1]{BD}.

For a kernel $\varphi $ on $\mathbb{R}$, write 
\begin{equation*}
B_{t}(\varphi ;f,g)(x)=\int_{\mathbb{R}}\varphi (u)f(x-tu)g(x+tu)\,du.
\end{equation*}%
For $a>0$ and $\tau \in \mathbb{R}$, set 
\begin{equation*}
D_{a}\varphi (u)=a^{-1}\varphi (u/a),\qquad T_{\tau }\varphi (u)=\varphi
(u-\tau ).
\end{equation*}%
With this notation, 
\begin{equation*}
A_{t}^{+}(f,g)=B_{t}(\mathbf{1}_{[0,1]};f,g),\qquad B_{t}(D_{a}\varphi
;f,g)=B_{at}(\varphi ;f,g).
\end{equation*}%
Let $\mathcal{S}_{0}$ denote the class of mean-zero functions $\theta $ such
that 
\begin{equation*}
\sup_{u\in \mathbb{R}}|u|^{a}\bigl|\frac{d_{{}}^{b}}{du^{b}}\theta (u)\bigr|%
\leq 1,\qquad 0\leq a,b\leq 100,
\end{equation*}%
where $a$ and $b$ are integers. For $\tau \in \mathbb{R}$, set 
\begin{equation*}
\mathcal{S}_{0}^{\tau }=\{T_{\tau }\theta :\theta \in \mathcal{S}_{0}\}.
\end{equation*}

\begin{thm}[{\protect\cite[Theorem~1.1]{BD}}]
\label{thm:BD-shifted} Let $1<r,s\le\infty$, $2/3<p<\infty$, and $%
1/p=1/r+1/s $. Let $m\ge1$ be an integer, and let $(\tau_j)_{j\in\mathbb{Z}}$
be a sequence of integers satisfying 
\begin{equation*}
2^{m-1}\le |\tau_j|\le 2^m,\qquad j\in\mathbb{Z}.
\end{equation*}
If $\psi_j\in\mathcal{S}_0^{\tau_j}$ for every $j\in\mathbb{Z}$, then for
all $f,g\in\mathcal{S}(\mathbb{R})$, 
\begin{equation*}
\Big\|\sum_{j\in\mathbb{Z}}B_{2^j}(\psi_j;f,g)\Big\|_{L^p(\mathbb{R})} \le
C_{p,r,s}\,m^4 \|f\|_{L^r(\mathbb{R})}\|g\|_{L^s(\mathbb{R})}.
\end{equation*}
\end{thm}

We first establish a maximal estimate for shifted bilinear averages with a
nonnegative smooth kernel.

\begin{lem}
\label{lem:shifted-maximal} Let $1<r,s<\infty$, $p>2/3$, and $1/p=1/r+1/s$.
Let $\Phi\in C_c^\infty(\mathbb{R})$ be nonnegative and set $%
\rho=\min\{1,p\} $. Let $h\ge0$ and $j_0$ be integers, and let $J\subseteq%
\mathbb{Z}$ be nonempty. For each $j\in J$, set $a_j=2^{j+j_0}$ and let $%
\tau_j\in[2^h,2^{h+1})$. Then, for all $f\in L^r(\mathbb{R})$ and $g\in L^s(%
\mathbb{R})$, 
\begin{equation*}
\Big\|\sup_{j\in J} B_{a_j}(T_{\tau_j}\Phi;|f|,|g|) \Big\|_{L^p(\mathbb{R}%
)}\le C_{\Phi,p,r,s}(1+h)^{4+1/\rho}\|f\|_{L^r(\mathbb{R})}\|g\|_{L^s(%
\mathbb{R})}.
\end{equation*}
The constant is independent of $h$, $J$, $j_0$, and the shifts $%
(\tau_j)_{j\in J}$.
\end{lem}

\begin{proof}
By approximation and monotone convergence, it suffices to assume that $f,g\in%
\mathcal{S}(\mathbb{R})$ and that $J$ is finite.

For the preliminary square-function estimate, suppose that $\tau _{j}\in
\lbrack 2^{h},2^{h+1})\cap \mathbb{Z}$ for every $j\in J$. Fix a compact
interval $K\subset \lbrack -R,R]$ independent of $j$. Let $(\theta
_{j})_{j\in J}$ be a family in $C_{c}^{\infty }(\mathbb{R})$ such that 
\begin{equation*}
\func{supp}\theta _{j}\subseteq K,\qquad \int_{\mathbb{R}}\theta
_{j}(u)\,du=0,\qquad j\in J.
\end{equation*}%
Assume also that 
\begin{equation*}
\sup_{j\in J}\Vert \theta _{j}^{(n)}\Vert _{L^{\infty }(\mathbb{R})}\leq
C_{n},\qquad n=0,\ldots ,100.
\end{equation*}%
The uniform support and derivative bounds give a constant $C_{\ast }>0$,
independent of $j$, such that $\theta _{j}/C_{\ast }\in \mathcal{S}_{0}$ for
every $j\in J.$ Fix a finite set $J\subset \mathbb{Z}$ and a choice of signs 
$\varepsilon _{j}\in \{-1,1\}$ for $j\in J$. Set 
\begin{equation*}
\psi _{j}(u)=%
\begin{cases}
C_{\ast }^{-1}\varepsilon _{j}\theta _{j}(u-\tau _{j}), & j\in J, \\ 
0, & j\notin J,%
\end{cases}%
\end{equation*}%
where $\tau _{j}$ denotes the integer shift. By construction, $\psi _{j}\in 
\mathcal{S}_{0}^{\tau _{j}}$ for every $j\in J$. Applying Theorem~\ref%
{thm:BD-shifted} with $m=h+1$ and reindexing the sum by $k=j+j_{0}$, we
obtain 
\begin{equation*}
\bigg\|\sum_{j\in J}B_{a_{j}}(\psi _{j};f,g)\bigg\|_{L^{p}(\mathbb{R}%
)}\lesssim _{p,r,s}(1+h)^{4}\Vert f\Vert _{L^{r}(\mathbb{R})}\Vert g\Vert
_{L^{s}(\mathbb{R})},\qquad a_{j}=2^{j+j_{0}}.
\end{equation*}%
Taking the $\varepsilon _{j}$ to be independent Rademacher variables, we
obtain from Khintchine's inequality and Tonelli's theorem that 
\begin{align}
\bigg\|\Big(\sum_{j\in J}\bigl|B_{a_{j}}(T_{\tau _{j}}\theta _{j};f,g)\bigr|%
^{2}\big)^{1/2}\bigg\|_{L^{p}(\mathbb{R})}& \lesssim _{p}C_{\ast }\bigg(%
\mathbb{E}_{\varepsilon }\Big\|\sum_{j\in J}B_{a_{j}}(\psi _{j};f,g)\Big\|%
_{L^{p}(\mathbb{R})}^{p}\bigg)^{1/p}  \notag \\
& \lesssim _{p,r,s}C_{\ast }(1+h)^{4}\Vert f\Vert _{L^{r}(\mathbb{R})}\Vert
g\Vert _{L^{s}(\mathbb{R})}.  \label{eq2}
\end{align}%
We now return to the real shifts in the statement of the lemma. Set $\Theta
=\Phi -D_{2}\Phi $, so that $\int_{\mathbb{R}}\Theta =0$. We then have 
\begin{equation*}
B_{a_{j}}(T_{\tau _{j}}\Phi ;|f|,|g|)=B_{a_{j}}(T_{\tau _{j}}\Theta
;|f|,|g|)+B_{2a_{j}}(T_{\tau _{j}/2}\Phi ;|f|,|g|).
\end{equation*}%
Applying the same identity with $a_{j}$ and $\tau _{j}$ replaced by $%
2^{n}a_{j}$ and $\tau _{j}/2^{n}$, respectively, we obtain 
\begin{equation*}
\begin{aligned} B_{2^na_j}(T_{\tau_j/2^n}\Phi;|f|,|g|)
-B_{2^{n+1}a_j}(T_{\tau_j/2^{n+1}}\Phi;|f|,|g|) =
B_{2^na_j}(T_{\tau_j/2^n}\Theta;|f|,|g|). \end{aligned}
\end{equation*}%
Summing over $n=0,\ldots ,h-1$ gives 
\begin{equation*}
B_{a_{j}}(T_{\tau _{j}}\Phi ;|f|,|g|)=\sum_{n=0}^{h-1}B_{2^{n}a_{j}}(T_{\tau
_{j}/2^{n}}\Theta ;|f|,|g|)+B_{2^{h}a_{j}}(T_{\tau _{j}/2^{h}}\Phi ;|f|,|g|).
\end{equation*}%
When $h=0$, the sum is understood to be empty. Fix $0\leq n<h$ and write 
\begin{equation*}
\frac{\tau _{j}}{2^{n}}=d_{j,n}+\delta _{j,n},\qquad d_{j,n}=\left\lfloor 
\frac{\tau _{j}}{2^{n}}\right\rfloor ,\qquad 0\leq \delta _{j,n}<1.
\end{equation*}%
Set $\Theta _{j,n}(u)=\Theta (u-\delta _{j,n})$. Then 
\begin{equation*}
T_{\tau _{j}/2^{n}}\Theta =T_{d_{j,n}}\Theta _{j,n},\qquad 2^{h-n}\leq
d_{j,n}<2^{h-n+1}.
\end{equation*}%
The kernels $\Theta _{j,n}$ have mean zero, common compact support, and
uniform derivative bounds. Applying (\ref{eq2}) at the scales $%
2^{n}a_{j}=2^{j+j_{0}+n}$, we obtain 
\begin{equation*}
\begin{aligned} \Big\|\sup_{j\in
J}\big|B_{2^na_j}(T_{\tau_j/2^n}\Theta;|f|,|g|)\big|\Big\|_{L^p} &\le
\bigg\|\Big(\sum_{j\in J}
\bigl|B_{2^na_j}(T_{d_{j,n}}\Theta_{j,n};|f|,|g|)\bigr|^2
\Big)^{1/2}\bigg\|_{L^p}\\
&\lesssim_{\Phi,p,r,s}(1+h-n)^4\|f\|_{L^r}\|g\|_{L^s}. \end{aligned}
\end{equation*}%
For the last term, note that $\tau _{j}/2^{h}\in \lbrack 1,2)$. If $\func{%
supp}\Phi \subset \lbrack -R,R]$, then 
\begin{equation*}
\func{supp}(T_{\tau _{j}/2^{h}}\Phi )\subset \lbrack -R-2,R+2].
\end{equation*}%
Changing variables by $y=2^{h}a_{j}u$, we obtain 
\begin{equation*}
\begin{aligned} B_{2^ha_j}(T_{\tau_j/2^h}\Phi;|f|,|g|)(x)
&\le\frac{\|\Phi\|_{L^\infty}}{2^ha_j} \int_{|y|\le
(R+2)2^ha_j}|f(x-y)g(x+y)|\,dy\\ &\le2(R+2)\|\Phi\|_{L^\infty}\,\mathcal
M(f,g)(x), \end{aligned}
\end{equation*}%
where 
\begin{equation*}
\mathcal{M}(f,g)(x)=\sup_{a>0}\frac{1}{2a}\int_{-a}^{a}|f(x-y)g(x+y)|\,dy.
\end{equation*}%
Taking the supremum over $j\in J$ and using Lacey's bilinear maximal
estimate \cite{Lacey}, we obtain 
\begin{equation*}
\Big\|\sup_{j\in J}B_{2^{h}a_{j}}(T_{\tau _{j}/2^{h}}\Phi ;|f|,|g|)\Big\|%
_{L^{p}(\mathbb{R})}\lesssim _{\Phi ,p,r,s}\Vert f\Vert _{L^{r}(\mathbb{R}%
)}\Vert g\Vert _{L^{s}(\mathbb{R})}.
\end{equation*}%
Combining the preceding estimates and using the $\rho $-subadditivity of the 
$L^{p}$ quasi-norm, where $\rho =\min \{1,p\}$, we obtain 
\begin{equation*}
\begin{aligned} \Big\|\sup_{j\in J}
B_{a_j}(T_{\tau_j}\Phi;|f|,|g|)\Big\|_{L^p(\mathbb R)}^\rho
&\lesssim_{\Phi,p,r,s}
\Big(1+\sum_{n=0}^{h-1}(1+h-n)^{4\rho}\Big)\|f\|_{L^r(\mathbb R)}^\rho
\|g\|_{L^s(\mathbb R)}^\rho\\ &\lesssim_{\Phi,p,r,s}
(1+h)^{4\rho+1}\|f\|_{L^r(\mathbb R)}^\rho\|g\|_{L^s(\mathbb R)}^\rho.
\end{aligned}
\end{equation*}%
Taking the $1/\rho $ power gives the desired estimate. The extension to
arbitrary $J\subseteq \mathbb{Z}$ follows by monotone convergence.
\end{proof}

This following counting lemma gives a bound for the $\ell ^{q}$ norm of a
finite-size array. It uses the largest entry of the array together with
averages obtained by picking either a single entry or a dyadic block from
each row. Its proof relies on a dyadic stopping-time argument.

\begin{lem}
\label{lem:array-counting} Let $p>0$ and $q>2$. Let $N=2^{m}$ for some
integer $m\geq 0$, and let $J\subseteq \mathbb{Z}$ be a nonempty finite set.
Consider a complex array $(a_{j,k})_{j\in J,\,0\leq k<N}$. For each $j\in J$%
, the entries $a_{j,0},\ldots ,a_{j,N-1}$ form the $j$ th row.

For each $M\in\{1,2,4,\ldots,N\}$, partition the column indices into the
dyadic blocks 
\begin{equation*}
I_{M,\ell}=\{\ell M,\ell M+1,\ldots,(\ell+1)M-1\},\qquad 0\le\ell<N/M.
\end{equation*}
Suppose that for each row $j$ and each dyadic block $I$, there is a
nonnegative number $P_{j,I}$ such that 
\begin{equation*}
\frac1{|I|}\sum_{k\in I}|a_{j,k}|\le C_0P_{j,I},
\end{equation*}
where $C_0>0$ is independent of $j$ and $I$.

In each row $j$, choose a column index $K_j$ uniformly from $%
\{0,\ldots,N-1\} $, with the choices independent across rows. For each $M$,
let $I_j^{(M)}$ be the unique dyadic block of length $M$ containing $K_j$.
For fixed $M$, each block of length $M$ is selected with probability $M/N$
in each row, independently across rows. Set 
\begin{equation*}
b=\max_{\substack{ j\in J  \\ 0\le k<N}}|a_{j,k}|, \qquad W=\bigg(\mathbb{E}%
\Big[\big(\sum_{j\in J}|a_{j,K_j}|^2\big)^{p/2}\Big]\bigg)^{1/p},
\end{equation*}
and 
\begin{equation*}
U=(m+1)^{1/p}\max_{M\in\{1,2,4,\ldots,N\}} \Big(\mathbb{E}\Big[\big(%
\max_{j\in J}P_{j,I_j^{(M)}}\big)^p\Big]\Big)^{1/p}.
\end{equation*}
Then 
\begin{equation*}
\Big(\sum_{j\in J}\sum_{k=0}^{N-1}|a_{j,k}|^q\Big)^{1/q}
\lesssim_{p,q,C_0}N^{1/q}\big(b^{1-1/q}U^{1/q}+b^{1-2/q}W^{2/q}\big).
\end{equation*}
\end{lem}

\begin{proof}
The conclusion is immediate if $b=0$. Assume that $b>0$ and fix $\lambda>0$.
For each $j\in J$, set 
\begin{equation*}
A_j(\lambda)=\{k\in\{0,\ldots,N-1\}:|a_{j,k}|>\lambda\}, \qquad
K_\lambda=\sum_{j\in J}|A_j(\lambda)|.
\end{equation*}
Thus $K_\lambda$ is the total number of entries in the array whose modulus
exceeds $\lambda$.

Suppose first that $0<K_\lambda\le N$. In each row $j$, we perform a dyadic
stopping-time decomposition by selecting all dyadic blocks $I$ that are
maximal with respect to inclusion among those satisfying 
\begin{equation*}
\frac{|A_j(\lambda)\cap I|}{|I|}\ge \frac{K_\lambda}{2N}.
\end{equation*}
If $I$ has a dyadic parent $\hat{I}$, then by the maximality of $I$ and the
identity $|\hat{I}|=2|I|$, we have 
\begin{equation*}
\frac{|A_j(\lambda)\cap I|}{|I|}\le2\frac{|A_j(\lambda)\cap \hat{I}|}{|\hat{I%
}|}<\frac{K_\lambda}{N}.
\end{equation*}
If $I$ is the entire row, then $I=\{0,\ldots,N-1\}$ and $|I|=N$. Thus, 
\begin{equation*}
\frac{|A_j(\lambda)\cap I|}{|I|} =\frac{|A_j(\lambda)|}{N}\le\frac{K_\lambda%
}{N}.
\end{equation*}
Let $\mathcal{I}_j$ be the family of selected blocks in row $j$. These
blocks are pairwise disjoint. Every singleton $\{k\}$ with $k\in
A_j(\lambda) $ satisfies the stopping condition, so the selected blocks
cover $A_j(\lambda)$. Hence 
\begin{equation*}
\begin{aligned} K_\lambda=\sum_{j\in J}\sum_{I\in\mathcal I_j}
|A_j(\lambda)\cap I|\le \frac{K_\lambda}{N} \sum_{j\in J}\sum_{I\in\mathcal
I_j}|I|. \end{aligned}
\end{equation*}
Dividing by $K_\lambda>0$, we obtain 
\begin{equation*}
\sum_{j\in J}\sum_{I\in\mathcal{I}_j}|I|\ge N.
\end{equation*}
For each $M\in\{1,2,4,\ldots,N\}$, let 
\begin{equation*}
C_M:=\sum_{j\in J}\#\{I\in\mathcal{I}_j:|I|=M\}
\end{equation*}
be the total number of selected blocks of length $M$ across all rows. Then 
\begin{equation*}
\sum_{M\in\{1,2,4,\ldots,N\}} M C_M=\sum_{j\in J}\sum_{I\in\mathcal{I}%
_j}|I|\ge N.
\end{equation*}
By the pigeonhole principle, there is a block length $M$ for which at least $%
N/[M(m+1)]$ blocks of length $M$ are selected across all rows. Let $%
c_j=\#\{I\in \mathcal{I}_j:|I|=M\}$. Then 
\begin{equation*}
\sum_{j\in J}c_j=C_M\ge\frac{N}{M(m+1)}.
\end{equation*}
In row $j$, the block $I_j^{(M)}$ is chosen uniformly from the $N/M$ blocks
of length $M$. Thus 
\begin{equation*}
\mathbb{P}\bigl(I_j^{(M)}\in\mathcal{I}_j\bigr)=\frac{c_j}{N/M}=\frac{Mc_j}{N%
}.
\end{equation*}
Let 
\begin{equation*}
X=\sum_{j\in J}\mathbf{1}_{\{I_j^{(M)}\in\mathcal{I}_j\}}
\end{equation*}
be the number of rows in which a retained block is sampled. Then 
\begin{equation*}
\mathbb{E }X = \sum_{j\in J}\mathbb{P}\bigl(I_j^{(M)}\in\mathcal{I}_j\bigr) %
= \frac{M}{N}\sum_{j\in J}c_j = \frac{MC_M}{N} \ge \frac1{m+1}.
\end{equation*}
For each $j\in J$, write $X_j=\mathbf{1}_{\{I_j^{(M)}\in\mathcal{I}_j\}}$,
so that $X=\sum_{j\in J}X_j$. The variables $X_j$ are independent across
rows and satisfy $X_j^2=X_j$. Hence 
\begin{equation*}
\func{Var}(X)=\sum_{j\in J}\func{Var}(X_j) =\sum_{j\in J}\left(\mathbb{E }%
X_j-(\mathbb{E }X_j)^2\right)\le \mathbb{E }X.
\end{equation*}
It follows that $\mathbb{E }X^2\le(\mathbb{E }X)^2+\mathbb{E }X$. Using $X=X%
\mathbf{1}_{\{X>0\}}$ and the Cauchy--Schwarz inequality, we obtain 
\begin{equation*}
\mathbb{E }X\le (\mathbb{E }X^2)^{1/2}\mathbb{P}(X>0)^{1/2}.
\end{equation*}
As $\mathbb{E }X\ge 1/(m+1)$, the probability of sampling at least one
retained block satisfies 
\begin{equation*}
\mathbb{P}(X>0)\ge \frac{(\mathbb{E }X)^2}{\mathbb{E }X^2} \ge \frac{\mathbb{%
E }X}{1+\mathbb{E }X}\ge \frac1{m+2}.
\end{equation*}
On the event $\{X>0\}$, choose a row $j_*\in J$ such that $I_{j_*}^{(M)}\in%
\mathcal{I}_{j_*}$. By the block average bound and the stopping condition,
we have 
\begin{equation*}
\begin{aligned} \max_{j\in J}P_{j,I_j^{(M)}}\ge P_{j_*,I_{j_*}^{(M)}}\ge
\frac{1}{C_0M}\sum_{k\in I_{j_*}^{(M)}}|a_{j_*,k}|\ge \frac{\lambda}{C_0M}
\bigl|A_{j_*}(\lambda)\cap I_{j_*}^{(M)}\bigr|\ge \frac{\lambda
K_\lambda}{2C_0N}. \end{aligned}
\end{equation*}
Set $Z_M=\max_{j\in J}P_{j,I_j^{(M)}}$. On the event $\{X>0\}$, we have $%
Z_M\ge \lambda K_\lambda/(2C_0N)$. Therefore, 
\begin{equation*}
\begin{aligned} \mathbb E[Z_M^p]\ge \mathbb
E\bigl[Z_M^p\mathbf1_{\{X>0\}}\bigr]\ge \left(\frac{\lambda
K_\lambda}{2C_0N}\right)^p \mathbb P(X>0)\ge\frac{1}{m+2}
\left(\frac{\lambda K_\lambda}{2C_0N}\right)^p. \end{aligned}
\end{equation*}
Taking the $p$th root and using the definition of $U$, we obtain 
\begin{equation*}
\begin{aligned} \lambda K_\lambda\le 2C_0N(m+2)^{1/p} \bigl(\mathbb
E[Z_M^p]\bigr)^{1/p}\le 2C_0N
\left(\frac{m+2}{m+1}\right)^{1/p}U\lesssim_{p,C_0} NU. \end{aligned}
\end{equation*}

Now suppose that $K_\lambda>N$. In each row $j\in J$, choose a column index $%
K_j$ independently and uniformly from $\{0,\ldots,N-1\}$. Let 
\begin{equation*}
X_\lambda=\sum_{j\in J}\mathbf{1}_{\{K_j\in A_j(\lambda)\}}
\end{equation*}
be the number of sampled entries whose modulus exceeds $\lambda$. Then 
\begin{equation*}
\mathbb{E }X_\lambda=\sum_{j\in J}\frac{|A_j(\lambda)|}{N}=\frac{K_\lambda}{N%
}>1.
\end{equation*}
The variable $X_\lambda$ is a sum of independent $\{0,1\}$-valued random
variables. As above, $\func{Var}(X_\lambda) \le \mathbb{E }X_\lambda$. Using 
$\mathbb{E }X_\lambda>1$, we obtain 
\begin{equation*}
\mathbb{E}[X_\lambda^2]=(\mathbb{E }X_\lambda)^2+\func{Var}(X_\lambda)\le (%
\mathbb{E }X_\lambda)^2+\mathbb{E }X_\lambda\le 2(\mathbb{E }X_\lambda)^2.
\end{equation*}
Set $E=\left\{X_\lambda\ge\mathbb{E }X_\lambda/2\right\}$. On $E^c$, we have 
$X_\lambda<\mathbb{E }X_\lambda/2$. Therefore, 
\begin{equation*}
\begin{aligned} \mathbb E X_\lambda=\mathbb E[X_\lambda\mathbf1_{E^c}]
+\mathbb E[X_\lambda\mathbf1_E]\le \frac12\mathbb E X_\lambda +\bigl(\mathbb
E[X_\lambda^2]\bigr)^{1/2}\mathbb P(E)^{1/2}, \end{aligned}
\end{equation*}
where the second term is bounded by the Cauchy--Schwarz inequality.
Subtracting $\frac12\mathbb{E }X_\lambda$ from both sides and squaring, we
obtain 
\begin{equation*}
\mathbb{P}(E)\ge \frac{(\mathbb{E }X_\lambda)^2} {4\mathbb{E}[X_\lambda^2]}%
\ge \frac18,
\end{equation*}
where the last step uses $\mathbb{E}[X_\lambda^2]\le 2(\mathbb{E }%
X_\lambda)^2$. On $E$, each sampled entry counted by $X_\lambda$ contributes
at least $\lambda^2$ to the sum of squares. Hence 
\begin{equation*}
\sum_{j\in J}|a_{j,K_j}|^2\ge \lambda^2 X_\lambda \ge \frac{\lambda^2}{2}%
\mathbb{E }X_\lambda=\frac{\lambda^2 K_\lambda}{2N}.
\end{equation*}
By the definition of $W$ and the bound $\mathbb{P}(E)\ge 1/8$, we get 
\begin{equation*}
\begin{aligned} W^p=\mathbb E\Big[ \big(\sum_{j\in
J}|a_{j,K_j}|^2\big)^{p/2} \Big]\ge\left(\frac{\lambda^2
K_\lambda}{2N}\right)^{p/2} \mathbb P(E)\ge \frac18 \left(\frac{\lambda^2
K_\lambda}{2N}\right)^{p/2}. \end{aligned}
\end{equation*}
Taking the $p$th root and rearranging, we obtain 
\begin{equation*}
\lambda^2 K_\lambda\lesssim_p N W^2.
\end{equation*}
By the distribution function formula, 
\begin{equation*}
\sum_{j\in J}\sum_{k=0}^{N-1}|a_{j,k}|^q =q\int_0^b
\lambda^{q-1}K_\lambda\,d\lambda.
\end{equation*}
On the set where $0<K_\lambda\le N$, we use $\lambda
K_\lambda\lesssim_{p,C_0}NU$. On the set where $K_\lambda>N$, we use $%
\lambda^2K_\lambda\lesssim_p NW^2$. The integrand vanishes when $K_\lambda=0$%
. Hence 
\begin{equation*}
\begin{aligned} \sum_{j\in J}\sum_{k=0}^{N-1}|a_{j,k}|^q &\lesssim_{p,q,C_0}
N\left( U\int_0^b\lambda^{q-2}\,d\lambda +W^2\int_0^b\lambda^{q-3}\,d\lambda
\right)\\ &\lesssim_{p,q,C_0} N\left(b^{q-1}U+b^{q-2}W^2\right).
\end{aligned}
\end{equation*}
The second integral converges because $q>2$. Taking the $q$th root, we
obtain 
\begin{equation*}
\left(\sum_{j\in J}\sum_{k=0}^{N-1}|a_{j,k}|^q\right)^{1/q}
\lesssim_{p,q,C_0} N^{1/q}\left(b^{1-1/q}U^{1/q}+b^{1-2/q}W^{2/q}\right),
\end{equation*}
as desired.
\end{proof}

Combining Lemma~\ref{lem:array-counting} with the maximal and
square-function estimates above, we obtain the following $L^p(\ell^q)$
estimate for an array of shifted bilinear averages.

\begin{lem}
\label{lem:shifted-array} Let $1<r,s<\infty$, $p>2/3$, and $1/p=1/r+1/s$.
Let $q>2$ and $N=2^m$, where $m\ge0$ is an integer. Let $J\subset\mathbb{Z}$
be a nonempty finite set and let $j_0\in\mathbb{Z}$.

Suppose that $\theta_k\in C_c^\infty(\mathbb{R})$ for $0\le k<N$ and that
there is a fixed compact interval $K$ such that $\func{supp}%
\theta_k\subseteq K$ and $\int_{\mathbb{R}}\theta_k(u)\,du=0.$ Assume also
that the bounds 
\begin{equation*}
\sup_{0\le k<N}\|\theta_k^{(n)}\|_{L^\infty(\mathbb{R})} \le C_n,\qquad
n=0,\ldots,100,
\end{equation*}
are uniform in $N$. Set 
\begin{equation*}
a_j=2^{j+j_0}, \qquad F_{j,k}=B_{a_j}(T_{N+k}\theta_k;f,g), \qquad j\in
J,\quad 0\le k<N.
\end{equation*}
Choose $\gamma$ such that 
\begin{equation*}
\max\left\{\frac1r,\frac1s,\frac{2}{3p}\right\}<\gamma<1, \qquad
\sigma=\frac1q+\left(1-\frac1q\right)\gamma.
\end{equation*}
Then $\sigma<1$ and 
\begin{equation*}
\Bigl\| \Bigl(\sum_{j\in J}\sum_{k=0}^{N-1}|F_{j,k}|^q\Bigr)^{1/q} \Bigr\|%
_{L^p(\mathbb{R})} \lesssim (1+m)^C N^\sigma \|f\|_{L^r(\mathbb{R}%
)}\|g\|_{L^s(\mathbb{R})}.
\end{equation*}
The constants are independent of $N$, $J$, and $j_0$, and their dependence
on the kernels is through $K$ and the bounds $C_n$.
\end{lem}

\begin{proof}
For each fixed $x\in\mathbb{R}$, apply Lemma~\ref{lem:array-counting} with $%
a_{j,k}=F_{j,k}(x)$. We estimate the corresponding quantities $b$, $W$, and $%
U$ separately.

We first estimate $b$. Let $K\subset[-R,R]$ be the fixed compact interval
containing the support of every $\theta_k$. Then 
\begin{equation*}
\func{supp}(T_{N+k}\theta_k)\subset K+(N+k)\subset[N+k-R,N+k+R]\subset[-CN,CN%
],
\end{equation*}
where $C$ is independent of $N$ and $k$. The uniform $L^\infty$ bound for $%
\theta_k$ and the fixed length of $K$ give 
\begin{equation*}
\|T_{N+k}\theta_k\|_{L^1(\mathbb{R})}+\|T_{N+k}\theta_k\|_{L^\infty(\mathbb{R%
})}\lesssim 1.
\end{equation*}
We arrive at 
\begin{equation*}
\begin{aligned} |F_{j,k}(x)|&= \Bigl|\int_{\mathbb
R}T_{N+k}\theta_k(u)f(x-a_ju)g(x+a_ju)\,du\Bigr|\\ &\le\int_{\mathbb
R}|T_{N+k}\theta_k(u)||f(x-a_ju)g(x+a_ju)|\,du\\
&\le\|T_{N+k}\theta_k\|_{L^1}^{1-\gamma} \Bigl(\int_{\mathbb
R}|T_{N+k}\theta_k(u)||f(x-a_ju)g(x+a_ju)|^{1/\gamma}\,du\Bigr)^\gamma\\
&\lesssim \Bigl(\frac1{a_j}\int_{|y|\le
CNa_j}|f(x-y)g(x+y)|^{1/\gamma}\,dy\Bigr)^\gamma\\ &\lesssim
N^\gamma\mathcal M\bigl(|f|^{1/\gamma},|g|^{1/\gamma}\bigr)(x)^\gamma.
\end{aligned}
\end{equation*}
The first inequality follows from the triangle inequality for integrals, and
the second follows from H\"older's inequality. Taking the maximum over $j\in
J$ and $0\le k<N$, we obtain the corresponding pointwise bound for $b$. Our
choice of $\gamma$ ensures that $\gamma r,\gamma s>1$, $\gamma p>2/3$, and 
\begin{equation*}
\frac1{\gamma p}=\frac1{\gamma r}+\frac1{\gamma s}.
\end{equation*}
Thus, by the bilinear maximal estimate \cite{Lacey}, 
\begin{equation*}
\begin{aligned} \|b\|_{L^p}&\lesssim N^\gamma \bigl\|\mathcal
M\bigl(|f|^{1/\gamma},|g|^{1/\gamma}\bigr) \bigr\|_{L^{\gamma
p}}^\gamma\lesssim N^\gamma \bigl\||f|^{1/\gamma}\bigr\|_{L^{\gamma
r}}^\gamma \bigl\||g|^{1/\gamma}\bigr\|_{L^{\gamma
s}}^\gamma=N^\gamma\|f\|_{L^r}\|g\|_{L^s}. \end{aligned}
\end{equation*}

Next, we estimate $W$. Fix a choice of column indices $(K_j)_{j\in J}$, one
in each row. The selected operator in row $j$ is 
\begin{equation*}
F_{j,K_j}=B_{a_j}(T_{N+K_j}\theta_{K_j};f,g).
\end{equation*}
The shifts $N+K_j$ are integers in $[N,2N)$. Applying the (\ref{eq2}) with $%
h=m$, we obtain 
\begin{equation*}
\Bigl\|\Bigl(\sum_{j\in J}|F_{j,K_j}|^2\Bigr)^{1/2}\Bigr\|_{L^p(\mathbb{R})}
\lesssim(1+m)^4\|f\|_{L^r(\mathbb{R})}\|g\|_{L^s(\mathbb{R})},
\end{equation*}
uniformly in the choice of $(K_j)_{j\in J}$. Averaging the $p$th power of
this estimate over all choices of $(K_j)_{j\in J}$ and using the definition
of $W$, we obtain 
\begin{equation*}
\begin{aligned} \|W\|_{L^p(\mathbb R)}^p=\mathbb E\bigg[
\Big\|\Big(\sum_{j\in J}|F_{j,K_j}|^2\Big)^{1/2}\Big\|_{L^p(\mathbb
R)}^p\bigg] \lesssim(1+m)^{4p}\|f\|_{L^r(\mathbb R)}^p\|g\|_{L^s(\mathbb
R)}^p. \end{aligned}
\end{equation*}
Taking the $p$th root gives 
\begin{equation*}
\|W\|_{L^p(\mathbb{R})}\lesssim (1+m)^4\|f\|_{L^r(\mathbb{R})}\|g\|_{L^s(%
\mathbb{R})}.
\end{equation*}

Finally, we estimate $U$. Fix a nonnegative function $\Phi\in C_c^\infty(%
\mathbb{R})$ such that $\Phi(u)>1$ for all $u\in[-R,R+1]$. Let $I$ be a
dyadic block of length $M$ starting at $k_0$, i.e., $I=\{k_0,k_0+1,%
\ldots,k_0+M-1\}$, where $k_0$ is a multiple of $M$ and $0\le k_0\le N-M$.
For each $k\in I$, the support of the shifted kernel satisfies 
\begin{equation*}
\func{supp}(T_{N+k}\theta_k)\subset [N+k-R,N+k+R]\subset
[N+k_0-R,N+k_0+M-1+R].
\end{equation*}
Fix $u\in\mathbb{R}$. If $u\in\func{supp}(T_{N+k}\theta_k)$, then 
\begin{equation*}
|u-(N+k)|\le R,\qquad\text{so}\qquad u-N-R\le k\le u-N+R.
\end{equation*}
This interval contains at most $\lfloor 2R\rfloor+1$ integers. Thus, at each 
$u$, only a bounded number of shifted kernels can be nonzero, uniformly in $%
N $ and $M$. For each $k\in I$, the support condition gives 
\begin{equation*}
|T_{N+k}\theta_k(u)|\le \|\theta_k\|_{L^\infty}\mathbf{1}_{[N+k-R,%
\,N+k+R]}(u).
\end{equation*}
Summing over $k\in I$, we obtain 
\begin{equation*}
\begin{aligned} \sum_{k\in I}|T_{N+k}\theta_k(u)|&\le\Bigl(\sup_{0\le
k<N}\|\theta_k\|_{L^\infty}\Bigr)\sum_{k\in
I}\mathbf1_{[N+k-R,\,N+k+R]}(u)\\ &\le\Bigl(\sup_{0\le
k<N}\|\theta_k\|_{L^\infty}\Bigr)(\lfloor
2R\rfloor+1)\mathbf1_{[N+k_0-R,\,N+k_0+M-1+R]}(u)\\
&\lesssim\Phi\Bigl(\frac{u-N-k_0}{M}\Bigr). \end{aligned}
\end{equation*}
The second inequality uses the bound $\lfloor 2R\rfloor+1$ on the number of
intervals containing $u$. For the last inequality, the normalized variable $%
(u-N-k_0)/M$ lies in $[-R,R+1]$ whenever the indicator is nonzero, and $%
\Phi>1$ on this interval. It follows that 
\begin{equation*}
\begin{aligned} \frac1M\sum_{k\in I}|F_{j,k}(x)|&\le\frac1M\int_{\mathbb
R}\Bigl(\sum_{k\in I}|T_{N+k}\theta_k(u)|\Bigr)|f(x-a_ju)g(x+a_ju)|\,du\\
&\lesssim\frac1M\int_{\mathbb
R}\Phi\Bigl(\frac{u-N-k_0}{M}\Bigr)|f(x-a_ju)g(x+a_ju)|\,du\\
&=B_{Ma_j}\bigl(T_{(N+k_0)/M}\Phi;|f|,|g|\bigr)(x)=:P_{j,I}(x), \end{aligned}
\end{equation*}
where the last equality follows from the change of variables $u=Mv$.

Fix a block length $M$ and a choice of blocks $(I_j^{(M)})_{j\in J}$. Let $%
k_{0,j}$ be the first column index in $I_j^{(M)}$. We have 
\begin{equation*}
Ma_j=2^{j+j_0+\log_2 M},\qquad\frac{N+k_{0,j}}{M}\in\Big[\frac NM,\frac{2N}{M%
}\Big).
\end{equation*}
Applying Lemma~\ref{lem:shifted-maximal} with $h=\log_2(N/M)\le m$, we
obtain 
\begin{equation*}
\Bigl\|\max_{j\in J}P_{j,I_j^{(M)}}\Bigr\|_{L^p(\mathbb{R})}
\lesssim_{\Phi,p,r,s}(1+m)^{4+1/\rho}\|f\|_{L^r(\mathbb{R})}\|g\|_{L^s(%
\mathbb{R})},
\end{equation*}
uniformly in $M$ and in the choice of $(I_j^{(M)})_{j\in J}$. For each fixed 
$M$, we take the $p$th power of the preceding estimate and average over all
choices of the blocks. By the definition of $U$, we obtain 
\begin{equation*}
\begin{aligned} \|U\|_{L^p(\mathbb R)}^p&\le
(m+1)\sum_{M\in\{1,2,4,\ldots,N\}} \mathbb E\Bigl[\Bigl\|\max_{j\in
J}P_{j,I_j^{(M)}}\Bigr\|_{L^p(\mathbb R)}^p\Bigr]\\
&\lesssim(1+m)^{p(4+1/\rho)+2}\|f\|_{L^r(\mathbb R)}^p\|g\|_{L^s(\mathbb
R)}^p, \end{aligned}
\end{equation*}
where the sum runs over the $m+1$ dyadic block lengths. Taking the $p$th
root gives 
\begin{equation*}
\|U\|_{L^p(\mathbb{R})}\lesssim(1+m)^{4+1/\rho+2/p}\|f\|_{L^r(\mathbb{R}%
)}\|g\|_{L^s(\mathbb{R})}.
\end{equation*}
Applying Lemma~\ref{lem:array-counting} and using H\"older's inequality in
the integral of the $p$th power of each product, we obtain 
\begin{equation*}
\begin{aligned} \Bigl\|\Bigl(\sum_{j\in
J}\sum_{k=0}^{N-1}|F_{j,k}|^q\Bigr)^{1/q}\Bigr\|_{L^p} \lesssim
N^{1/q}\Bigl(\|b\|_{L^p}^{1-1/q}\|U\|_{L^p}^{1/q}+\|b\|_{L^p}^{1-2/q}\|W%
\|_{L^p}^{2/q} \Bigr). \end{aligned}
\end{equation*}
Here we use $p$-subadditivity to combine the two terms when $p<1$.
Substituting the preceding bounds for $b$, $U$, and $W$, we find that the
powers of $N$ in the two terms are 
\begin{equation*}
\frac1q+\Bigl(1-\frac1q\Bigr)\gamma=\sigma,\qquad\frac1q+\Bigl(1-\frac2q%
\Bigr)\gamma=\sigma-\frac{\gamma}{q},
\end{equation*}
respectively. Thus 
\begin{equation*}
\Bigl\|\Bigl(\sum_{j\in J}\sum_{k=0}^{N-1}|F_{j,k}|^q\Bigr)^{1/q}\Bigr\|%
_{L^p}\lesssim (1+m)^C N^\sigma\|f\|_{L^r}\|g\|_{L^s},
\end{equation*}
as desired.
\end{proof}

The next two lemmas establish short-variation estimates for a fixed smooth
kernel and for a family of shifted smooth kernels, respectively.

\begin{lem}
\label{lem:smooth-short} Let $1<r,s<\infty$, $p>2/3$, and $1/p=1/r+1/s$.
Then, for every $q>2$ and every $\varphi\in C_c^\infty(\mathbb{R})$, we have 
\begin{equation*}
\Big\|V^{S,q}\bigl(B_t(\varphi;f,g)\bigr)\Big\|_{L^p(\mathbb{R})}
\lesssim_{\varphi,p,q,r,s}\|f\|_{L^r(\mathbb{R})}\|g\|_{L^s(\mathbb{R})}.
\end{equation*}
\end{lem}

\begin{proof}
We may assume that $f,g\in\mathcal{S}(\mathbb{R})$. For each $j\in\mathbb{Z}$%
, write $t=2^j c$ with $1\le c\le 2$ on the dyadic interval $[2^j,2^{j+1}]$.
The dilation identity gives $B_{2^j c}(\varphi;f,g)=B_{2^j}(D_c\varphi;f,g).$
For $1\le c\le 2$, the kernels $D_c\varphi$ and $\partial_cD_c\varphi$ are
uniformly bounded and supported in a common compact interval.
Differentiating under the integral gives 
\begin{equation*}
\frac{d}{dc}B_{2^j}(D_c\varphi;f,g)(x)=B_{2^j}(\partial_cD_c\varphi;f,g)(x).
\end{equation*}
Fix a finite increasing sequence $2^j\le t_0<\cdots<t_N\le 2^{j+1}$ and set $%
c_k=2^{-j}t_k$. By the dilation identity and the fundamental theorem of
calculus, 
\begin{equation*}
B_{t_{k+1}}(\varphi;f,g)(x)-B_{t_k}(\varphi;f,g)(x)
=\int_{c_k}^{c_{k+1}}B_{2^j}(\partial_cD_c\varphi;f,g)(x)\,dc.
\end{equation*}
Taking absolute values and using $\|\cdot\|_{\ell^q}\le\|\cdot\|_{\ell^1}$,
we obtain 
\begin{equation*}
\begin{aligned} \Bigl(\sum_{k=0}^{N-1} \bigl|B_{t_{k+1}}(\varphi;f,g)(x)
-B_{t_k}(\varphi;f,g)(x)\bigr|^q\Bigr)^{1/q}
&\le\sum_{k=0}^{N-1}\int_{c_k}^{c_{k+1}}
\bigl|B_{2^j}(\partial_cD_c\varphi;f,g)(x)\bigr|\,dc\\
&\le\int_1^2\bigl|B_{2^j}(\partial_cD_c\varphi;f,g)(x)\bigr|\,dc.
\end{aligned}
\end{equation*}
Here the intervals $[c_k,c_{k+1}]$ partition $[c_0,c_N]\subseteq[1,2]$.
Taking the supremum over all such sequences yields 
\begin{equation*}
V_q\bigl(B_t(\varphi;f,g):2^j\le t\le 2^{j+1}\bigr)(x) \le\int_1^2\bigl|%
B_{2^j}(\partial_cD_c\varphi;f,g)(x)\bigr|\,dc.
\end{equation*}
The dilation $D_c\varphi(u)=c^{-1}\varphi(u/c)$ preserves the integral, so $%
\int_{\mathbb{R}}\partial_cD_c\varphi(u)\,du=0$. Choose $\zeta\in
C_c^\infty((1/2,5/2))$ with $\zeta=1$ on $[1,2]$. For $c>0$ and $u\in\mathbb{%
R}$, set 
\begin{equation*}
G(c,u)=\zeta(c)\,\partial_cD_c\varphi(u).
\end{equation*}
Set $G(0,u)=0$. By the support condition on $\zeta$, $G(c,u)=0$ for $c\in[%
0,1/2]\cup[5/2,3]$. Thus the restriction of $G$ to $[0,3]\times\mathbb{R}$
extends smoothly and $3$-periodically in $c$. For $\ell\in\mathbb{Z}$, set 
\begin{equation*}
h_\ell(u)=\frac13\int_0^3 G(c,u)e^{-2\pi i\ell c/3}\,dc.
\end{equation*}
The Fourier expansion of $G$ in $c$ gives 
\begin{equation*}
G(c,u)=\sum_{\ell\in\mathbb{Z}}h_\ell(u)e^{2\pi i\ell c/3}.
\end{equation*}
By Fubini's theorem and the mean-zero property of $\partial_cD_c\varphi$, we
have 
\begin{equation*}
\int_{\mathbb{R}}h_\ell(u)\,du=\frac13\int_0^3 e^{-2\pi i\ell c/3} \Bigl(%
\int_{\mathbb{R}}G(c,u)\,du\Bigr)\,dc=0,\qquad \ell\in\mathbb{Z}.
\end{equation*}
Choose $R>0$ such that $\func{supp}\varphi\subseteq[-R,R]$. For $c\in\func{%
supp}\zeta$, the kernel $\partial_cD_c\varphi$ is supported in $[-cR,cR]%
\subseteq[-5R/2,5R/2]$. It follows that 
\begin{equation*}
\func{supp}h_\ell\subseteq[-5R/2,5R/2],\qquad \ell\in\mathbb{Z}.
\end{equation*}
Let $A\ge1$ and $n\ge0$ be integers. For $\ell\ne0$, integrating by parts $A$
times in $c$ gives 
\begin{equation*}
\partial_u^n h_\ell(u)=\frac13\Bigl(\frac{3}{2\pi i\ell}\Bigr)^A\int_0^3
\partial_c^A\partial_u^n G(c,u)e^{-2\pi i\ell c/3}\,dc.
\end{equation*}
All boundary terms vanish because $G$ vanishes near $c=0$ and $c=3$. Thus $%
\|\partial_u^n h_\ell\|_{L^\infty}\lesssim_{\varphi,A,n}|\ell|^{-A}$ for $%
\ell\ne0$. For $\ell=0$, we estimate directly 
\begin{equation*}
\|\partial_u^n h_0\|_{L^\infty}\le \frac13\int_0^3\|\partial_u^n
G(c,\cdot)\|_{L^\infty}\,dc\lesssim_{\varphi,n}1.
\end{equation*}
Combining the two cases, we obtain 
\begin{equation*}
\|\partial_u^n
h_\ell\|_{L^\infty}\lesssim_{\varphi,A,n}(1+|\ell|)^{-A},\qquad \ell\in%
\mathbb{Z}.
\end{equation*}
On $[1,2]$, we have $\zeta(c)=1$, so that $G(c,u)=\partial_cD_c\varphi(u)$.
The common compact support and rapid decay of the coefficients $h_\ell$
allow us to substitute the Fourier series into $B_{2^j}$ term by term. Thus 
\begin{equation*}
\begin{aligned} V_q\bigl(B_t(\varphi;f,g):2^j\le t\le 2^{j+1}\bigr)(x)
&\le\int_1^2\big|\sum_{\ell\in\mathbb Z}e^{2\pi i\ell c/3}
B_{2^j}(h_\ell;f,g)(x)\big|\,dc\\ &\le\sum_{\ell\in\mathbb
Z}|B_{2^j}(h_\ell;f,g)(x)|. \end{aligned}
\end{equation*}
By the triangle inequality in $\ell^q$ and the inequality $%
\|\cdot\|_{\ell^q}\le\|\cdot\|_{\ell^2}$ for $q>2$, we obtain 
\begin{equation*}
\begin{aligned} V^{S,q}\bigl(B_t(\varphi;f,g)\bigr)(x) &\le
\sum_{\ell\in\mathbb Z} \Bigl(\sum_{j\in\mathbb
Z}|B_{2^j}(h_\ell;f,g)(x)|^q\Bigr)^{1/q}\\ &\le \sum_{\ell\in\mathbb Z}
\Bigl(\sum_{j\in\mathbb Z}|B_{2^j}(h_\ell;f,g)(x)|^2\Bigr)^{1/2}.
\end{aligned}
\end{equation*}
Fix $\ell\in\mathbb{Z}$ and set $\widetilde h_\ell(u)=h_\ell(u+1)$. In the
square-function estimate \eqref{eq2}, take 
\begin{equation*}
\theta_j=(1+|\ell|)^A\widetilde h_\ell,\qquad \tau_j=1,\qquad a_j=2^j
\end{equation*}
for every $j$. The kernels $\theta_j$ have mean zero, a common compact
support, and derivative bounds independent of $\ell$. Applying \eqref{eq2}
with $h=0$ to finite sets of indices $j$ and then using monotone
convergence, we obtain 
\begin{equation*}
\Bigl\|\Bigl(\sum_{j\in\mathbb{Z}}|B_{2^j}(h_\ell;f,g)|^2\Bigr)^{1/2}\Bigr\|%
_{L^p(\mathbb{R})} \lesssim_{\varphi,A,p,r,s}(1+|\ell|)^{-A}\|f\|_{L^r(%
\mathbb{R})}\|g\|_{L^s(\mathbb{R})}.
\end{equation*}
Set $\rho=\min\{1,p\}$ and choose an integer $A$ such that $A\rho>1$. By the 
$\rho$-subadditivity of the $L^p$ quasi-norm and the preceding
square-function estimate, we have 
\begin{equation*}
\begin{aligned}
\Bigl\|V^{S,q}\bigl(B_t(\varphi;f,g)\bigr)\Bigr\|_{L^p}^{\rho}&\le
\sum_{\ell\in\mathbb Z} \Bigl\|\Bigl(\sum_{j\in\mathbb
Z}|B_{2^j}(h_\ell;f,g)|^2\Bigr)^{1/2}\Bigr\|_{L^p}^{\rho}\\
&\lesssim_{\varphi,A,p,r,s}\|f\|_{L^r}^{\rho}\|g\|_{L^s}^{\rho}\sum_{\ell\in%
\mathbb Z}(1+|\ell|)^{-A\rho}. \end{aligned}
\end{equation*}
The series converges because $A\rho>1$. Taking the $1/\rho$ power gives the
desired estimate.
\end{proof}

\begin{lem}
\label{lem:shifted-short} Let $1<r,s<\infty$, $p>2/3$, and $1/p=1/r+1/s$.
Let $q>2$ and $\psi\in C_c^\infty(\mathbb{R})$ satisfy $\int_{\mathbb{R}%
}\psi(u)\,du=0$. For $N=2^m$, where $m\ge1$ is an integer, we have 
\begin{equation*}
\Big\|V^{S,q}\bigl(B_{t/N}(T_N\psi;f,g)\bigr)\Big\|_{L^p(\mathbb{R})}
\lesssim(1+m)^C N^\sigma\|f\|_{L^r(\mathbb{R})}\|g\|_{L^s(\mathbb{R})}.
\end{equation*}
Here $\sigma<1$ is the exponent from Lemma~\ref{lem:shifted-array}, and both 
$C$ and the implicit constant are independent of $m$.
\end{lem}

\begin{proof}
Fix $j\in\mathbb{Z}$ and divide the dyadic interval $[2^j,2^{j+1}]$ into $N$
subintervals of equal length. Denote the $k$th subinterval by 
\begin{equation*}
I_{j,k}=\Bigl[2^j\Bigl(1+\frac{k}{N}\Bigr),\, 2^j\Bigl(1+\frac{k+1}{N}\Bigr)%
\Bigr],\qquad k=0,\ldots,N-1.
\end{equation*}
For $t\in I_{j,k}$, write 
\begin{equation*}
t=2^j c,\qquad c=1+\frac{k+w}{N},\qquad 0\le w\le1.
\end{equation*}
Write $u=N+k+v$. The identity $Nc=N+k+w$ gives 
\begin{equation*}
\frac{u}{c}-N=\frac{u-Nc}{c}=\frac{v-w}{c}.
\end{equation*}
Consequently, 
\begin{equation*}
(D_cT_N\psi)(u)=(D_cT_N\psi)(N+k+v)=\frac1c\psi\Bigl(\frac{v-w}{c}\Bigr)%
=(T_wD_c\psi)(v).
\end{equation*}
Choose $\zeta\in C_c^\infty((-1/4,5/4))$ with $\zeta=1$ on $[0,1]$. For $w\in%
[-1/2,3/2]$ and $v\in\mathbb{R}$, set $c=1+(k+w)/N$ and define 
\begin{equation*}
H_{k,N}(w,v)=\zeta(w)(T_wD_c\psi)(v)=\frac{\zeta(w)}{c}\,\psi\Bigl(\frac{v-w%
}{c}\Bigr).
\end{equation*}
By the support condition on $\zeta$, $H_{k,N}$ vanishes near both endpoints
of $[-1/2,3/2]$. We therefore extend $H_{k,N}$ smoothly and $2$-periodically
in $w$. Let $\theta_{k,\ell}$ denote the $\ell$th Fourier coefficient of $%
H_{k,N}$ in $w$, given by 
\begin{equation*}
\theta_{k,\ell}(v)=\frac12\int_{-1/2}^{3/2}H_{k,N}(w,v)e^{-\pi i\ell
w}\,dw,\qquad \ell\in\mathbb{Z}.
\end{equation*}
The change of variables $z=(v-w)/c$ and Fubini's theorem give 
\begin{equation*}
\int_{\mathbb{R}}\theta_{k,\ell}(v)\,dv=0.
\end{equation*}
For every $A>0$ and every integer $n\ge0$, repeated integration by parts in $%
w$ for $\ell\ne0$ and a direct estimate for $\ell=0$ yield 
\begin{equation*}
\|\partial_v^n\theta_{k,\ell}\|_{L^\infty(\mathbb{R})}\lesssim_{%
\psi,A,n}(1+|\ell|)^{-A},
\end{equation*}
with an implicit constant independent of $N$, $k$, and $\ell$. For $t\in
I_{j,k}$, we have $w\in[0,1]$ and hence $\zeta(w)=1$. By the dilation
identity and the preceding translation formula, 
\begin{equation*}
B_{t/N}(T_N\psi;f,g)=B_{2^{j-m}}(D_cT_N\psi;f,g)=B_{2^{j-m}}\bigl(%
T_{N+k}H_{k,N}(w,\cdot);f,g\bigr).
\end{equation*}
Substituting the Fourier series of $H_{k,N}$ term by term, we obtain 
\begin{equation*}
B_{t/N}(T_N\psi;f,g)=\sum_{\ell\in\mathbb{Z}}e^{\pi i\ell
w}B_{2^{j-m}}(T_{N+k}\theta_{k,\ell};f,g).
\end{equation*}
Fix $x\in\mathbb{R}$. In the preceding Fourier expansion, the coefficients $%
B_{2^{j-m}}(T_{N+k}\theta_{k,\ell};f,g)(x)$ are independent of $w$. Using $%
|e^{\pi i\ell w}|=1$, we obtain 
\begin{equation*}
\sup_{t\in I_{j,k}}\bigl|B_{t/N}(T_N\psi;f,g)(x)\bigr|\le\sum_{\ell\in%
\mathbb{Z}} \bigl|B_{2^{j-m}}(T_{N+k}\theta_{k,\ell};f,g)(x)\bigr|.
\end{equation*}
Take any finite increasing sequence $t_0<\cdots<t_M$ in $I_{j,k}$ and write 
\begin{equation*}
t_\nu=2^j\Bigl(1+\frac{k+w_\nu}{N}\Bigr),\qquad 0\le w_0<\cdots<w_M\le1.
\end{equation*}
By the mean value inequality, for each $\ell\in\mathbb{Z}$, 
\begin{equation*}
\sum_{\nu=0}^{M-1}\bigl|e^{\pi i\ell w_{\nu+1}}-e^{\pi i\ell w_\nu}\bigr|\le
\pi|\ell|(w_M-w_0)\le \pi|\ell|.
\end{equation*}
Subtracting the Fourier expansions at $t_{\nu+1}$ and $t_\nu$, taking
absolute values, and summing over $\nu$, we obtain a bound uniform over all
such sequences. Taking the supremum yields 
\begin{equation*}
V_1\bigl(B_{t/N}(T_N\psi;f,g):t\in I_{j,k}\bigr)(x)\le \pi\sum_{\ell\in%
\mathbb{Z}}|\ell|\,\bigl|B_{2^{j-m}}(T_{N+k}\theta_{k,\ell};f,g)(x)\bigr|.
\end{equation*}
Combining the two estimates gives 
\begin{equation*}
\begin{aligned} &\sup_{t\in
I_{j,k}}\bigl|B_{t/N}(T_N\psi;f,g)(x)\bigr|+V_1\bigl(B_{t/N}(T_N\psi;f,g):t%
\in I_{j,k}\bigr)(x)\\ &\qquad\lesssim\sum_{\ell\in\mathbb
Z}(1+|\ell|)\,\bigl|B_{2^{j-m}}(T_{N+k}\theta_{k,\ell};f,g)(x)\bigr|.
\end{aligned}
\end{equation*}
Fix $j\in\mathbb{Z}$ and $x\in\mathbb{R}$, and take any finite increasing
sequence 
\begin{equation*}
2^j\le t_0<\cdots<t_M\le 2^{j+1}.
\end{equation*}
Let $k_\nu$ be the index of the subinterval containing $t_\nu$. Assign each
shared endpoint to the subinterval on its right, and assign $2^{j+1}$ to the
last subinterval. Thus the indices are uniquely determined and satisfy $%
k_0\le\cdots\le k_M$.

For a fixed $k$, the partition points assigned to $I_{j,k}$ form a
consecutive subsequence. The sum of the absolute increments between
successive points of this subsequence is bounded by the total variation on $%
I_{j,k}$. Using $\|\cdot\|_{\ell^q}\le\|\cdot\|_{\ell^1}$, we obtain 
\begin{equation*}
\begin{aligned} &\sum_{\substack{0\le\nu<M\\ k_\nu=k_{\nu+1}=k}}
\Bigl|B_{t_{\nu+1}/N}(T_N\psi;f,g)(x)-B_{t_\nu/N}(T_N\psi;f,g)(x)\Bigr|^q\\
&\qquad\le V_1\bigl(B_{t/N}(T_N\psi;f,g):t\in I_{j,k}\bigr)(x)^q.
\end{aligned}
\end{equation*}
For increments whose endpoints lie in different subintervals, we use $%
|a-b|^q\le 2^{q-1}(|a|^q+|b|^q)$ to obtain 
\begin{equation*}
\begin{aligned} &\sum_{\substack{0\le\nu<M\\ k_\nu<k_{\nu+1}}}
\Bigl|B_{t_{\nu+1}/N}(T_N\psi;f,g)(x)-B_{t_\nu/N}(T_N\psi;f,g)(x)\Bigr|^q\\
&\qquad\le 2^{q-1}\sum_{\substack{0\le\nu<M\\ k_\nu<k_{\nu+1}}}
\sum_{a=\nu}^{\nu+1}\bigl|B_{t_a/N}(T_N\psi;f,g)(x)\bigr|^q\le
2^q\sum_{k=0}^{N-1} \sup_{t\in
I_{j,k}}\bigl|B_{t/N}(T_N\psi;f,g)(x)\bigr|^q. \end{aligned}
\end{equation*}
Summing the first estimate over $k$, adding the second estimate, and taking
the supremum over all such sequences, we obtain 
\begin{equation*}
\begin{aligned} &V_q\bigl(B_{t/N}(T_N\psi;f,g):2^j\le t\le
2^{j+1}\bigr)(x)^q\\
&\qquad\le\sum_{k=0}^{N-1}V_1\bigl(B_{t/N}(T_N\psi;f,g):t\in
I_{j,k}\bigr)(x)^q+ 2^q\sum_{k=0}^{N-1}\sup_{t\in
I_{j,k}}\bigl|B_{t/N}(T_N\psi;f,g)(x)\bigr|^q. \end{aligned}
\end{equation*}
Substituting the estimates on $I_{j,k}$ into the bound on each dyadic
interval and applying the triangle inequality in $\ell^q$ over $(j,k)$, we
obtain 
\begin{equation*}
\begin{aligned} &V^{S,q}\bigl(B_{t/N}(T_N\psi;f,g)\bigr)(x)\\
&\qquad\lesssim_q\sum_{\ell\in\mathbb Z}(1+|\ell|)\Bigl(\sum_{j\in\mathbb
Z}\sum_{k=0}^{N-1}
\bigl|B_{2^{j-m}}(T_{N+k}\theta_{k,\ell};f,g)(x)\bigr|^q\Bigr)^{1/q}.
\end{aligned}
\end{equation*}
Set $\rho=\min\{1,p\}$ and choose $A$ such that $(A-1)\rho>1$. For each
fixed $\ell$, the kernels $(1+|\ell|)^A\theta_{k,\ell}$ have mean zero, a
common compact support, and derivative bounds uniform in $N$, $k$, and $\ell$%
. Applying Lemma~\ref{lem:shifted-array} to these kernels with $j_0=-m$,
first for a finite set of indices $j$ and then by monotone convergence,
gives 
\begin{equation*}
\begin{aligned} \Bigl\|\Bigl(\sum_{j\in\mathbb
Z}\sum_{k=0}^{N-1}\bigl|B_{2^{j-m}}(T_{N+k}\theta_{k,\ell};f,g)\bigr|^q%
\Bigr)^{1/q}\Bigr\|_{L^p} \lesssim(1+|\ell|)^{-A}(1+m)^C
N^\sigma\|f\|_{L^r}\|g\|_{L^s}. \end{aligned}
\end{equation*}
By the $\rho$-subadditivity of the $L^p$ quasi-norm, we obtain 
\begin{equation*}
\begin{aligned}
\Bigl\|V^{S,q}\bigl(B_{t/N}(T_N\psi;f,g)\bigr)\Bigr\|_{L^p}^{\rho}
\lesssim\Bigl((1+m)^C
N^\sigma\|f\|_{L^r}\|g\|_{L^s}\Bigr)^\rho\sum_{\ell\in\mathbb
Z}(1+|\ell|)^{-(A-1)\rho}. \end{aligned}
\end{equation*}
The series converges because $(A-1)\rho>1$. Taking the $1/\rho$ power gives
the desired estimate.
\end{proof}

We now prove the short-variation estimate using the preceding two lemmas.

\begin{proof}[Proof of Proposition~\protect\ref{Becher-short}]
To decompose the averaging kernel $\mathbf{1}_{[0,1]}$, choose a nonnegative
even function $\eta\in C_c^\infty(\mathbb{R})$ such that $\int_{\mathbb{R}%
}\eta(u)\,du=1$. Define 
\begin{equation*}
\psi(u)=\int_{-\infty}^{u}\bigl(\eta(v)-D_2\eta(v)\bigr)\,dv.
\end{equation*}
The integrand is smooth and compactly supported with integral zero, so $%
\psi\in C_c^\infty(\mathbb{R})$. Its evenness gives $\psi(0)=0$, and hence $%
\psi$ is odd. In particular, $\int_{\mathbb{R}}\psi(u)\,du=0$. For $m\ge1$,
set $N=2^m$. By the definition of dilation, 
\begin{equation*}
(D_{1/N}\eta-D_{2/N}\eta)(y)=N(\eta-D_2\eta)(Ny).
\end{equation*}
Making the change of variables $v=N(u-y)$, we obtain 
\begin{equation*}
\begin{aligned} \bigl[\mathbf{1}_{[0,1]}*(D_{1/N}\eta-D_{2/N}\eta)\bigr](u)
&=N\int_0^1(\eta-D_2\eta)\bigl(N(u-y)\bigr)\,dy\\
&=\int_{N(u-1)}^{Nu}(\eta-D_2\eta)(v)\,dv\\ &=\psi(Nu)-\psi(Nu-N)\\
&=\frac1N\bigl[D_{1/N}(\psi-T_N\psi)\bigr](u). \end{aligned}
\end{equation*}
Summing over $m=1,\ldots,L$ and telescoping gives 
\begin{equation*}
\mathbf{1}_{[0,1]}*D_{2^{-L}}\eta=\mathbf{1}_{[0,1]}*\eta+%
\sum_{m=1}^{L}2^{-m}D_{2^{-m}}\psi-\sum_{m=1}^{L}2^{-m}D_{2^{-m}}(T_{2^m}%
\psi).
\end{equation*}
By linearity in the kernel and the dilation identity $B_t(D_a%
\varphi;f,g)=B_{at}(\varphi;f,g)$, we obtain 
\begin{equation*}
\begin{aligned}
B_t(\mathbf{1}_{[0,1]}*D_{2^{-L}}\eta;f,g)=&B_t(\mathbf{1}_{[0,1]}*%
\eta;f,g)+\sum_{m=1}^{L}2^{-m}B_{t/2^m}(\psi;f,g)\\
&-\sum_{m=1}^{L}2^{-m}B_{t/2^m}(T_{2^m}\psi;f,g). \end{aligned}
\end{equation*}
Since $\mathbf{1}_{[0,1]}*\eta\in C_c^\infty(\mathbb{R})$, Lemma~\ref%
{lem:smooth-short} gives 
\begin{equation*}
\Bigl\|V^{S,q}\bigl(B_t(\mathbf{1}_{[0,1]}*\eta;f,g)\bigr)\Bigr\|%
_{L^p}\lesssim \|f\|_{L^r}\|g\|_{L^s}.
\end{equation*}
The change of parameters $t\mapsto t/2^m$ maps $[2^j,2^{j+1}]$ onto $%
[2^{j-m},2^{j-m+1}]$. Reindexing the dyadic intervals and applying Lemma~\ref%
{lem:smooth-short} to $\psi$, we obtain 
\begin{equation*}
\begin{aligned}
\Bigl\|V^{S,q}\bigl(B_{t/2^m}(\psi;f,g)\bigr)\Bigr\|_{L^p}=\Bigl\|V^{S,q}%
\bigl(B_t(\psi;f,g)\bigr)\Bigr\|_{L^p}\lesssim \|f\|_{L^r}\|g\|_{L^s},
\end{aligned}
\end{equation*}
with an implicit constant independent of $m$. For the last term, the
mean-zero property of $\psi$ allows us to apply Lemma~\ref{lem:shifted-short}
with $N=2^m$, which gives 
\begin{equation*}
\Bigl\|V^{S,q}\bigl(B_{t/2^m}(T_{2^m}\psi;f,g)\bigr)\Bigr\|_{L^p}
\lesssim(1+m)^C2^{m\sigma}\|f\|_{L^r}\|g\|_{L^s}.
\end{equation*}
Set $\rho=\min\{1,p\}$. Applying the pointwise triangle inequality for the
short variation to the finite decomposition and using the $\rho$%
-subadditivity of the $L^p$ quasi-norm, we obtain 
\begin{equation*}
\begin{aligned}
&\Bigl\|V^{S,q}\bigl(B_t(\mathbf{1}_{[0,1]}*D_{2^{-L}}\eta;f,g)\bigr)\Bigr%
\|_{L^p}^{\rho}\\
&\quad\lesssim\Bigl(1+\sum_{m=1}^{L}2^{-m\rho}+\sum_{m=1}^{L}(1+m)^{C%
\rho}2^{-m(1-\sigma)\rho}\Bigr)\|f\|_{L^r}^{\rho}\|g\|_{L^s}^{\rho}\\
&\quad\lesssim \|f\|_{L^r}^{\rho}\|g\|_{L^s}^{\rho}. \end{aligned}
\end{equation*}
The last inequality follows from $\sigma<1$, and the bound is uniform in $L$%
. Finally, as $L\to\infty$, dominated convergence gives $B_t(\mathbf{1}%
_{[0,1]}*D_{2^{-L}}\eta;f,g)\to A_t^+(f,g)$. The lower semicontinuity of the
short variation and Fatou's lemma then give the desired estimate.
\end{proof}

\section{Variational Estimates for the Error Term}

\label{set 3} In this section we treat the error term $E_{\varepsilon }$,
which arises from the difference between sharp cutoff and smooth cutoff. The
key idea is to represent the singular kernel $D$ as an integral over a
finite-variation signed measure, so that its variation estimate reduces to
the variational bounds for bilinear averages obtained in Section 3.

\begin{thm}
\label{thm error} Let $1<r,s<\infty$, $p>2/3$, and $1/p=1/r+1/s$. Then, for
every $q>2$ we have 
\begin{equation*}
\Big\|V_q\big(E_{\varepsilon}(f,g):\varepsilon>0\big)\Big\|_{L^p}\leq
C_{p,q,r,s} \|f\|_{L^r}\|g\|_{L^s}.
\end{equation*}
\end{thm}

Before proving Theorem \ref{thm error}, we first analyze the error kernel 
\begin{eqnarray*}
D(u) &:&=\frac{1_{\{|u|>1\}}-\chi (u)}{u},\qquad u\neq 0, \\
D(0) &:&=0.
\end{eqnarray*}

The following lemma holds.

\begin{lem}\label{lem 4.2}
For almost every $u\in\mathbb{R}$, we have
\begin{align}
D(u)
&=-\func{sgn}(u)1_{\{|u|<1\}}
+\func{sgn}(u)\int_1^2
\left(\chi^{\prime}(t)+\frac{1-\chi(t)}{t}\right)
\frac{1_{\{|u|<t\}}}{t}\,dt.
\label{eq1}
\end{align}
Furthermore, there exists a finite signed Borel measure $\nu$ on
$(0,\infty)$, supported in $[1,2]$, with
\[
\|\nu\|_{\mathrm{TV}}<\infty,
\]
such that, for almost every $u\in\mathbb{R}$,
\[
D(u)=\int_{[1,2]}\frac{\func{sgn}(u)}{t}
1_{\{|u|<t\}}\,d\nu(t).
\]
\end{lem}

\begin{proof}
Since both sides of \eqref{eq1} are odd functions of $u$, it suffices
to consider $u\in[0,\infty)$. Note that
\[
\frac1t\left(\chi^{\prime}(t)+\frac{1-\chi(t)}t\right)
=-\frac{d}{dt}\left(\frac{1-\chi(t)}t\right).
\]
Using $\chi(1)=0$ and $\chi(2)=1$, for $1\leq a\leq2$, we obtain
\[
\int_a^2\left(\chi^{\prime}(t)+\frac{1-\chi(t)}t\right)
\frac{dt}{t}=\frac{1-\chi(a)}a.
\]
By our conventions, both sides vanish at $u=0$. When $0<u<1$,
we have $\chi(u)=0$ and hence $D(u)=0$. The right-hand side of
\eqref{eq1} equals
\[
-1+\int_1^2\left(\chi^{\prime}(t)+\frac{1-\chi(t)}t\right)
\frac{dt}{t}=-1+1=0.
\]
When $1<u<2$, the first term on the right-hand side of \eqref{eq1}
vanishes, while the integral term equals
\[
\begin{aligned}
\int_1^2\left(\chi'(t)+\frac{1-\chi(t)}t\right)
\frac{1_{\{u<t\}}}{t}\,dt
&=\int_u^2\left(\chi'(t)+\frac{1-\chi(t)}t\right)\frac{dt}{t}\\
&=\frac{1-\chi(u)}u=D(u).
\end{aligned}
\]
When $u\geq2$, we have $\chi(u)=1$ and hence $D(u)=0$.
The right-hand side of \eqref{eq1} also vanishes. Thus, the identity
holds on $\mathbb{R}\setminus\{-1,1\}$, and hence for almost every
$u\in\mathbb{R}$.

To obtain the measure representation, define
\[
d\nu(t):=-d\delta_1(t)+1_{(1,2)}(t)
\left(\chi^{\prime}(t)+\frac{1-\chi(t)}t\right)\,dt,
\]
where $\delta_1$ denotes the Dirac measure at $t=1$. This defines a
signed Borel measure $\nu$ supported in $[1,2]$, as $\delta_1$ is a
finite positive Borel measure and the density of the second term is
integrable. The atomic and absolutely continuous parts are mutually
singular. Thus, by the Jordan decomposition theorem,
\[
d|\nu|(t)=d\delta_1(t)+1_{(1,2)}(t)
\left|\chi^{\prime}(t)+\frac{1-\chi(t)}t\right|\,dt.
\]
Then
\[
\begin{aligned}
\|\nu\|_{\mathrm{TV}}
&=|\nu|((0,\infty))
=\int_{(0,\infty)}1\,d\delta_1(t)
+\int_{(0,\infty)}
\left|\chi'(t)+\frac{1-\chi(t)}{t}\right|
\mathbf{1}_{(1,2)}(t)\,dt\\
&=\delta_1((0,\infty))
+\int_1^2\left|\chi'(t)+\frac{1-\chi(t)}{t}\right|\,dt\\
&=1+\int_1^2\left|\chi'(t)+\frac{1-\chi(t)}{t}\right|\,dt\\
&\leq 1+\left\|\chi^\prime\right\|_{L^1(1,2)}
+\left\|1-\chi\right\|_{L^\infty(1,2)}\log 2<\infty.
\end{aligned}
\]
Thus, $\nu$ is finite, and the desired kernel representation follows
immediately from \eqref{eq1} and the definition of $\nu$.
\end{proof}

\begin{proof}[Proof of Theorem~\protect\ref{thm error}]
Using the definition of $D$, we can write the error term as 
\begin{equation*}
\begin{aligned}
E_{\varepsilon}(f,g)(x)&=\frac{1}{\varepsilon}\int_{%
\mathbb{R}}f(x-y)g(x+y)D(\frac{y}{\varepsilon})dy\\
&=\int_{\mathbb{R}}f(x-y)g(x+y)\int_{[1,2]}\frac{\operatorname{sgn}(y)}{t%
\varepsilon}1_{\{|y|<t\varepsilon\}}d\nu(t)dy. \end{aligned}
\end{equation*}
For $f,g\in\mathcal{S}(\mathbb{R})$, we have 
\begin{equation*}
\begin{aligned}
\int_{[1,2]}\frac1{t\varepsilon}\int_{|y|<t\varepsilon}|f(x-y)g(x+y)|\,dy%
\,d|\nu|(t)\le2\|f\|_{L^\infty}\|g\|_{L^\infty}\|\nu\|_{\mathrm{TV}}<\infty.
\end{aligned}
\end{equation*}
By Fubini's theorem, we obtain 
\begin{equation*}
\begin{aligned}
E_\varepsilon(f,g)(x)&=\int_{[1,2]}\frac1{t\varepsilon}\int_{|y|<t%
\varepsilon}\operatorname{sgn}(y)f(x-y)g(x+y)\,dy\,d\nu(t)\\ &=\int_{[1,2]}
\left[\frac1{t\varepsilon}\int_0^{t\varepsilon}f(x-y)g(x+y)\,dy-\frac1{t%
\varepsilon}\int_0^{t\varepsilon} f(x+y)g(x-y)\,dy \right]d\nu(t)\\
&=\int_{[1,2]}\bigl(A_{t\varepsilon}^+(f,g)(x)
-A_{t\varepsilon}^-(f,g)(x)\bigr)\,d\nu(t) \end{aligned}
\end{equation*}
the second equality follows from the change of variables $y\mapsto-y$. Fix
an arbitrary finite increasing sequence 
\begin{equation*}
0<\varepsilon_0<\cdots<\varepsilon_N.
\end{equation*}
The representation above yields 
\begin{equation*}
\begin{aligned}
E_{\varepsilon_{k+1}}(f,g)(x)-E_{\varepsilon_k}(f,g)(x)=&\int_{[1,2]}%
\Bigl[A^+_{t\varepsilon_{k+1}}(f,g)(x)-A^+_{t\varepsilon_k}(f,g)(x)\\
&\quad\quad\quad-A^-_{t\varepsilon_{k+1}}(f,g)(x)+A^-_{t%
\varepsilon_k}(f,g)(x)\Bigr]\,d\nu(t). \end{aligned}
\end{equation*}
Taking absolute values and then applying Minkowski's integral inequality in $%
\ell^q$, we obtain 
\begin{equation*}
\begin{aligned}
&\Bigl(\sum_{k=0}^{N-1}\Big|E_{\varepsilon_{k+1}}(f,g)(x)-E_{%
\varepsilon_k}(f,g)(x)\Big|^q\Bigr)^{\frac{1}{q}}\\ \leq
&\int_{[1,2]}\Bigl(\sum_{k=0}^{N-1}\Big|A^+_{t%
\varepsilon_{k+1}}(f,g)(x)-A^+_{t\varepsilon_k}(f,g)(x)-A^-_{t%
\varepsilon_{k+1}}(f,g)(x)+A^-_{t\varepsilon_k}(f,g)(x)\Big|^q\Bigr)^{%
\frac{1}{q}}\,d|\nu|(t)\\ \leq
&\int_{[1,2]}\Bigl[\Bigl(\sum_{k=0}^{N-1}\Big|A^+_{t%
\varepsilon_{k+1}}(f,g)(x)-A^+_{t\varepsilon_k}(f,g)(x)\Big|^q\Bigr)^{%
\frac{1}{q}}\\
&\qquad\qquad\qquad\qquad+\Bigl(\sum_{k=0}^{N-1}\Big|A^-_{t%
\varepsilon_{k+1}}(f,g)(x)-A^-_{t\varepsilon_k}(f,g)(x)\Big|^q\Bigr)^{%
\frac{1}{q}}\Bigr]\,d|\nu|(t). \end{aligned}
\end{equation*}
Now, we arrive at 
\begin{equation*}
\begin{aligned}
\left(\sum_{k=0}^{N-1}|E_{\varepsilon_{k+1}}-E_{\varepsilon_k}|^q
\right)^{1/q}\le \int_{[1,2]}
\Big[V_q(A_{t\varepsilon}^+:\varepsilon>0)+V_q(A_{t\varepsilon}^-:%
\varepsilon>0)\Big]\,d|\nu|(t). \end{aligned}
\end{equation*}
For each fixed $t>0$, the change of parameter $\varepsilon\mapsto
t\varepsilon$ gives 
\begin{equation*}
V_q\bigl(A_{t\varepsilon}^{\pm}:\varepsilon>0\bigr)=V_q\bigl(%
A_\varepsilon^{\pm}:\varepsilon>0\bigr).
\end{equation*}
Substituting this into the preceding inequality, we obtain 
\begin{equation*}
\begin{aligned} V_q\big(E_{\varepsilon}(f,g):\varepsilon>0\big)&\leq
\Bigl[V_q\big(A^+_{\varepsilon}(f,g):\varepsilon>0\big)+V_q\big(A^-_{%
\varepsilon}(f,g):\varepsilon>0\big)\Bigr]\int_{[1,2]}d|\nu|(t)\\
&=\|\nu\|_{\mathrm{TV}}\Bigl[V_q\big(A^+_{\varepsilon}(f,g):\varepsilon>0%
\big)+V_q\big(A^-_{\varepsilon}(f,g):\varepsilon>0\big)\Bigr]. \end{aligned}
\end{equation*}
Taking the $L^p$ quasi-norm and applying Theorem~\ref{Becher}, using $p$%
-subadditivity when $p<1$, completes the proof.
\end{proof}
\medskip
\noindent\textbf{Acknowledgements.}\enspace
This work is supported by the National Natural Science Foundation of
China through project numbers 12671125 and 12222114.

\end{document}